\documentclass{amsart}
\usepackage{amssymb}
\usepackage{graphicx}
\theoremstyle{plain}
\newtheorem{theorem}{Theorem}[section]
\newtheorem{corollary}[theorem]{Corollary}
\newtheorem{lemma}[theorem]{Lemma}
\newtheorem{proposition}[theorem]{Proposition}
\theoremstyle{definition}

\newtheorem{example}[theorem]{Example}
\newtheorem{conjecture}[theorem]{Conjecture}
\theoremstyle{remark}
\newtheorem{remark}[theorem]{Remark\upshape }

\newtheorem{acknowledgements}{Acknowledgements}

\numberwithin{equation}{section}
\renewcommand{\theenumi}{\roman{enumi}}

\newcommand{\at}{{@}}
\newcommand{\figlabone}{\relax}
\newcommand{\figlabtwo}{\relax}
\newlength{\qedskip}
\newlength{\displayboxwidth}
\newbox\ipbox
\newcommand{\ip}[2]{\left\langle #1\mathrel{\mathchoice
{\setbox\ipbox=\hbox{$\displaystyle \left\langle\mathstrut #1#2\right\rangle$}
\vrule height\ht\ipbox width0.25pt depth\dp\ipbox}
{\setbox\ipbox=\hbox{$\textstyle \left\langle\mathstrut #1#2\right\rangle$}
\vrule height\ht\ipbox width0.25pt depth\dp\ipbox}
{\setbox\ipbox=\hbox{$\scriptstyle \left\langle\mathstrut #1#2\right\rangle$}
\vrule height\ht\ipbox width0.25pt depth\dp\ipbox}
{\setbox\ipbox=\hbox{$\scriptscriptstyle \left\langle\mathstrut #1#2\right\rangle$}
\vrule height\ht\ipbox width0.25pt depth\dp\ipbox}
} #2\right\rangle}
\newcommand{\diracb}[1]{\left\langle #1\mathrel{\mathchoice
{\setbox\ipbox=\hbox{$\displaystyle \left\langle\mathstrut #1\right.$}
\vrule height\ht\ipbox width0.25pt depth\dp\ipbox}
{\setbox\ipbox=\hbox{$\textstyle \left\langle\mathstrut #1\right.$}
\vrule height\ht\ipbox width0.25pt depth\dp\ipbox}
{\setbox\ipbox=\hbox{$\scriptstyle \left\langle\mathstrut #1\right.$}
\vrule height\ht\ipbox width0.25pt depth\dp\ipbox}
{\setbox\ipbox=\hbox{$\scriptscriptstyle \left\langle\mathstrut #1\right.$}
\vrule height\ht\ipbox width0.25pt depth\dp\ipbox}
}\right. }
\newcommand{\dirack}[1]{\left. \mathrel{\mathchoice
{\setbox\ipbox=\hbox{$\displaystyle \left.\mathstrut #1\right\rangle$}
\vrule height\ht\ipbox width0.25pt depth\dp\ipbox}
{\setbox\ipbox=\hbox{$\textstyle \left.\mathstrut #1\right\rangle$}
\vrule height\ht\ipbox width0.25pt depth\dp\ipbox}
{\setbox\ipbox=\hbox{$\scriptstyle \left.\mathstrut #1\right\rangle$}
\vrule height\ht\ipbox width0.25pt depth\dp\ipbox}
{\setbox\ipbox=\hbox{$\scriptscriptstyle \left.\mathstrut #1\right\rangle$}
\vrule height\ht\ipbox width0.25pt depth\dp\ipbox}
} #1\right\rangle}
\newcommand{\rip}[2]{\left( #1\mathrel{\mathchoice
{\setbox\ipbox=\hbox{$\displaystyle \left(\mathstrut #1#2\right)$}
\vrule height\ht\ipbox width0.25pt depth\dp\ipbox}
{\setbox\ipbox=\hbox{$\textstyle \left(\mathstrut #1#2\right)$}
\vrule height\ht\ipbox width0.25pt depth\dp\ipbox}
{\setbox\ipbox=\hbox{$\scriptstyle \left(\mathstrut #1#2\right)$}
\vrule height\ht\ipbox width0.25pt depth\dp\ipbox}
{\setbox\ipbox=\hbox{$\scriptscriptstyle \left(\mathstrut #1#2\right)$}
\vrule height\ht\ipbox width0.25pt depth\dp\ipbox}
} #2\right)}

\def\limfunc#1{\mathop{\rm #1}}%
\long\def\TeXButton#1#2{#2}%
\def\LaTeXparent#1{}%
\def\ChildStyles#1{}%

\begin{document}
\title[Commuting self-adjoint symmetric partial derivatives]
{Commuting self-adjoint extensions of symmetric operators defined from the partial derivatives}
\author{Palle E. T. Jorgensen}
\address{Department of Mathematics\\
The University of Iowa\\
Iowa City, IA 52242\\
U.S.A.}
\email{jorgen%
\at%
math.uiowa.edu}
\author{Steen Pedersen}
\address{Department of Mathematics\\
Wright State University\\
Dayton, OH 45435\\
U.S.A.}
\email{steen%
\at%
math.wright.edu}
\thanks{Work supported by the National Science Foundation.}

\begin{abstract}
We consider the problem of
finding commuting self-adjoint
extensions of the partial derivatives
$\left\{ \frac{1}{i}\frac{\partial \;}{\partial x_j}:j=1,\dots ,d\right\} $ with domain
$C_{c}^{\infty}\left( \Omega \right) $ where the self-adjointness
is defined relative to $L^{2}\left( \Omega \right) $,
and $\Omega $ is a given open
subset of $\mathbb{R}^{d}$. The measure
on $\Omega $ is Lebesgue measure
on $\mathbb{R}^{d}$ restricted to $\Omega $.
The problem originates
with I.E. Segal and B.~Fuglede,
and is difficult in general.
In this paper, we provide
a representation-theoretic
answer in the special case
when $\Omega =I\times \Omega_{2}$ and
$I$ is an open interval. We
then apply the results to the
case when $\Omega $ is a $d$-cube, $I^{d}$,
and we describe
possible subsets $\Lambda \subset \mathbb{R}^{d}$
such that $\left\{ e_{\lambda}|_{I^{d}}:\lambda \in \Lambda \right\} $
is an orthonormal basis in
$L^{2}\left( I^{d}\right) $.
\end{abstract}

\maketitle\tableofcontents\listoffigures\renewcommand
{\figlabone}{\setlength{\unitlength}{54pt}\smash{\makebox[0pt]{\begin{picture}(4,5)(-2.71,0.02)
\put(-1,-0.1){\makebox(0.25,1)[lt]{$\beta _{1}-\beta _{0}$}}
\put(0,0.11){\makebox(0.25,1)[lt]{$\beta _{2}-\beta _{1}$}}
\put(1,0.36){\makebox(0.25,1)[lt]{$\beta _{3}-\beta _{2}$}}
\end{picture}}}}\renewcommand
{\figlabtwo}{\setlength{\unitlength}{54pt}\smash{\makebox[0pt]{\begin{picture}(5,4)(-0.65,-2.04)
\put(-0.1,-1){\makebox(0.25,1)[rt]{$\beta _{1}-\beta _{0}$}}
\put(0.11,0){\makebox(0.25,1)[rt]{$\beta _{2}-\beta _{1}$}}
\put(0.36,1){\makebox(0.25,1)[rt]{$\beta _{3}-\beta _{2}$}}
\end{picture}}}}\setlength{\displayboxwidth}{\textwidth}\addtolength{\displayboxwidth
}{-2\leftmargini}

\section{\label{S1}Introduction}

Recently several papers have appeared on
commuting non-self-adjoint operators
and their spectral theory; see, e.g.,
\cite{LKMV95}. The present paper
concerns the case when the given commuting operators
are unbounded and symmetric, but non-self-adjoint.
A concrete class of operators is studied, and
we address the questions of when commuting extension operators exist and,
when they do exist, what their structural properties are.

The problem of understanding
commuting symmetric, but
non-self-adjoint, unbounded
operators also has an origin in
mathematical physics
\cite{AvIv95,Bel89,HaKo91,Pav79}. The
terminology from physics is
``hermitian'', or \emph{``formally} self-adjoint'',
for symmetry, i.e., for the identity
$\ip{Sf}{h}=\ip{f}{Sh}$ for all vectors
$f,h$ in the domain of the operator $S$.
The simplest case of this is the
problem of assigning quantum mechanical boundary
conditions for free particles
confined in a box.
More
specifically, the problem here corresponds to the quantum-mechanical
trajectories of a particle confined to a region of tube type, e.g., a unit
cube. It is ``free'' except for the boundary conditions, and variations of
the boundary conditions (as considered here) correspond to different
physics.
For single
operators, von Neumann solved
(or made precise) the problem
by use
of the Cayley transform, and
considering instead the extension
problem for partial isometries.
But this approach does not
work well in the case of
several operators. Powers
(in \cite{Pow71,Pow74}) introduced an
algebraic approach for understanding
several operators, but the present
problem is very concrete and does
not lend itself easily to the
algebraic techniques introduced by
Powers.

Closely connected to the problem of finding commuting self-adjoint
extensions of $\frac{1}{i}\frac{\partial \;}{\partial x_{j}}$, $j=1,\dots ,d$%
, on $C_{c}^{\infty }\left( \Omega \right) $ is the corresponding spectral
question: If commuting self-adjoint extensions do exist, then it is known
that the common eigenfunctions of the extension operators must be of the
form $e_{\lambda }:=e^{i\lambda \cdot x}$ for special values of $\lambda \in 
\mathbb{R}^{d}$. Hence the spectral problem is that of finding when a given
pair $\left( \Omega ,\Lambda \right) $ satisfies the condition that $\left\{
e_{\lambda }|_{\Omega }:\lambda \in \Lambda \right\} $ is an orthogonal
basis in the Hilbert space $\mathcal{L}^{2}\left( \Omega \right) $. We note
that this so-called \emph{spectral pair} condition is very restrictive, and
so it explains the rigid geometric configurations $\left( \Omega ,\Lambda
\right) $ which admit solutions. But it also serves to motivate recent very
interesting developments on overcomplete systems; see, e.g., \cite%
{Kem99a,Kem99b}.

The setting of \emph{spectral pairs} in $d$ real dimensions involves two
subsets $\Omega $ and $\Lambda $ in $\mathbb{R}^{d}$ such that $\Omega $ has
finite and positive $d$-dimensional Lebesgue measure, and $\Lambda $ is an
index set for an orthogonal $\mathcal{L}^{2}\left( \Omega \right) $-basis $%
e_{\lambda }$ of exponentials, i.e., 
\begin{equation}
e_{\lambda }\left( x\right) =e^{i2\pi \lambda \cdot x},\quad x\in \Omega
,\;\lambda \in \Lambda   \label{eq1}
\end{equation}
where $\lambda \cdot x=\sum_{j=1}^{d}\lambda _{j}x_{j}$. We use vector
notation $x=\left( x_{1},\cdots ,x_{d}\right) $, $\lambda =\left( \lambda
_{1},\cdots ,\lambda _{d}\right) $, $x_{j},\lambda _{j}\in \mathbb{R}$, $%
j=1,\dots ,d$. The basis property refers to the Hilbert space $\mathcal{L}%
^{2}\left( \Omega \right) $ with inner product 
\begin{equation}
\ip{f}{g}%
_{\Omega }:=\int_{\Omega }\overline{f\left( x\right) }g\left( x\right) \,dx
\label{eq2}
\end{equation}
where $dx=dx_{1}\cdots dx_{d}$, and $f,g\in \mathcal{L}^{2}\left( \Omega
\right) $. The corresponding norm is 
\begin{equation}
\left\| f\right\| _{\Omega }^{2}:=%
\ip{f}{f}%
_{\Omega }=\int_{\Omega }\left| f\left( x\right) \right| ^{2}\,dx,
\label{eq3}
\end{equation}
as usual. It follows that the spectral pair property for a pair $\left(
\Omega ,\Lambda \right) $ is equivalent to 
\begin{equation*}
\Lambda -\Lambda =\left\{ \lambda -\lambda ^{\prime }:\lambda ,\lambda
^{\prime }\in \Lambda \right\} 
\end{equation*}
being contained in the \emph{zero-set} of the complex function 
\begin{equation}
z\longmapsto \int_{\Omega }e^{i2\pi z\cdot x}\,dx=:F_{\Omega }\left(
z\right)   \label{eq4}
\end{equation}
where $z=\left( z_{1},\cdots ,z_{d}\right) \in \mathbb{C}^{d}$, and $z\cdot
x:=\sum_{j=1}^{d}z_{j}x_{j}$, and the corresponding $e_{\lambda }$-set $%
\left\{ e_{\lambda }:\lambda \in \Lambda \right\} $ being \emph{total} in $%
\mathcal{L}^{2}\left( \Omega \right) $. Recall, totality means that the span
of the $e_{\lambda }$'s is dense in $\mathcal{L}^{2}\left( \Omega \right) $
relative to the $\left\| \,\cdot \,\right\| _{\Omega }$-norm, or,
equivalently, that $f=0$ is the only $\mathcal{L}^{2}\left( \Omega \right) $%
-solution to: 
\begin{equation*}
\ip{f}{e_\lambda }%
_{\Omega }=0\text{,\quad for all }\lambda \in \Lambda .
\end{equation*}

\section{Spectral pairs\label{SP}}

The theory of \emph{spectral pairs} 
was developed in 
previous joint papers by the
coauthors \cite{
JoPe92,
JoPe94,
JoPe96}.
A set $\Omega $ with finite nonzero Lebesgue measure is called a
\emph{spectral set} if $\left( \Omega ,\Lambda \right) $
is a spectral pair for some set $\Lambda $.
We recall that Fuglede showed
\cite{Fug74} that the disk and the triangle in
two dimensions are \emph{not spectral
sets.}
By the disk and the triangle
we mean the usual versions,
respectively, 
$\left\{ \left( x_{1}^{{}},x_{2}^{{}}\right) \in \mathbb{R}%
_{{}}^{2}:x_{1}^{2}+x_{2}^{2}<1\right\}  $
and
$\left\{ \left( x_{1},x_{2}\right) \in \mathbb{R}^{2}:0<x_{1},\;0<x_{2},%
\;x_{1}+x_{2}<1\right\} .$
\linebreak
Note that, for the present discussion, it is inessential whether or not the
sets $\Omega $ are taken to be open, but it is essential for the following
theorem which we will need. It is due to Fuglede and the coauthors; see \cite
{Fug74,Jor82,Ped87,JoPe92}.

If $\Omega \subset \mathbb{R}^{d}$ is open, then we consider the partial
derivatives $\frac{\partial \;}{\partial x_{j}}$, $j=1,\dots ,d$, defined on 
$C_{c}^{\infty }\left( \Omega \right) $ as unbounded skew-symmetric
operators in $\mathcal{L}^{2}\left( \Omega \right) $. The corresponding
versions $\frac{1}{2\pi \sqrt{-1}}\frac{\partial \;}{\partial x_{j}}$ are
symmetric of course. We say that $\Omega $ has the \emph{extension property}
if there are commuting self-adjoint extension operators $H_{j}$, i.e., 
\begin{equation}
\frac{1}{2\pi i}\frac{\partial \;}{\partial x_{j}}\subset H_{j},\quad j=1,\dots
,d.  \label{eq7}
\end{equation}
We say that the containment $A\subset B$ holds for two operators $A$ and $B$
if the graph of $A$ is contained in that of $B$. (For details, see \cite
{ReSi} and \cite{DS2}.) \emph{Commutativity} for the extension operators $%
H_{j}$ is in the strong sense of spectral resolutions.
Since the $H_{j}$'s are assumed self-adjoint, each one has a
projection-valued \emph{spectral resolution} $E_{j}$, i.e., an $\mathcal{L}%
^{2}\left( \Omega \right) $-projection-valued Borel measure on $\mathbb{R}$,
such that $E_{j}\left( \mathbb{R}\right) =I_{\mathcal{L}^{2}\left( \Omega
\right) }$, and 
\begin{equation}
H_{j}=\int_{-\infty }^{\infty }\lambda E_{j}\left( d\lambda \right) ,
\label{eq8}
\end{equation}
for $j=1,\dots ,d$. The strong commutativity is taken to mean 
\begin{equation}
E_{j}\left( \Delta \right) E_{j^{\prime }}\left( \Delta ^{\prime }\right)
=E_{j^{\prime }}\left( \Delta ^{\prime }\right) E_{j}\left( \Delta \right) 
\label{eq9}
\end{equation}
for all $j,j^{\prime }=1,\dots ,d$, and all Borel subsets $\Delta ,\Delta
^{\prime }\subset \mathbb{R}$.
Extensions commuting in a weaker sense were considered in \cite{Fri87}.

Our analysis is based on von Neumann's
deficiency-space characterization of the
self-adjoint extensions of a given symmetric
operator \cite{vNeu29}. Let $\Omega $ be an open set
with finite Lebesgue measure. For each $j$,
the deficiency spaces corresponding to
$\frac{1}{i}\frac{\partial \;}{\partial x_{j}}$
are infinite-dimensional. It follows
that each
$\frac{1}{i}\frac{\partial \;}{\partial x_{j}}$
has ``many'' self-adjoint
extensions. The main problem (not addressed
by von Neumann's theory) is the selection
of a commuting set
$H_{1},H_{2},\dots ,H_{d}$
of extensions.
In fact, for some $\Omega $ (e.g., when $d=2$,
the disk and the triangle) it is impossible
to select a commuting set
$H_{1},H_{2},\dots ,H_{d}$
of extensions.

We have (see \cite{Fug74,Jor82,Ped87,JoPe92})

\begin{theorem}
\label{Thm1.1}\textup{(Fuglede, Jorgensen, Pedersen)} Let $\Omega \subset %
\mathbb{R}^{d}$ be open and connected with finite and positive Lebesgue
measure. Then $\Omega $ has the extension property if and only if it is a
spectral set. Moreover, with $\Omega $ given, there is a one-to-one correspondence
between the two sets of subsets:
\begin{equation}
\left\{ \Lambda \subset {\mathbb{R}}^{d}\colon
\left( \Omega ,\Lambda \right) 
\text{ is a spectral pair}
\right\}
\label{eqThm1.1a}
\end{equation}
and
\begin{multline}
\{ \Lambda \subset {\mathbb{R}}^{d}\colon
\Lambda 
\text{ is the joint spectrum of some commutative}\\
\text{family }
\left( H_{1},\dots ,H_{d}\right) 
\text{ of self-adjoint etensions}
\}.
\label{eqThm1.1b}
\end{multline}
This correspondence is determined as follows:
\renewcommand{\theenumi}{\alph{enumi}}
\begin{enumerate}
\item \label{Thm1.1(1)}If the extensions
$\left( H_{1},\dots ,H_{d}\right) $ are
given, then $\lambda \in \Lambda $ if and only if
\begin{equation}
e_{\lambda}\in\bigcap_{j}\mathop{\rm domain}\left( H_{j}\right) .
\label{eqThm1.1(1)}
\end{equation}

\item \label{Thm1.1(2)}If,
conversely, $\left( \Omega ,\Lambda \right) $ is a spectral
pair at the outset, then the ansatz \textup{(\ref{eqThm1.1(1)})} and
\begin{equation}
H_{j}e_{\lambda}=\lambda_{j}e_{\lambda},\qquad \lambda \in \Lambda
\label{eqThm1.1(2)}
\end{equation}
determine uniquely a set of commuting extensions.
\end{enumerate}

If $\Omega $ is only assumed open, then the spectral-set
property implies the extension property, but not conversely.
\end{theorem}

\begin{corollary}
\label{RemNew1.2}
Suppose $\Omega$ is open and connected.
It follows then that
a discrete set $\Lambda $ is the joint
spectrum of some commuting self-adjoint extension
operators $H_{j}$, $j=1,\dots ,d$, if and only if 
$\left( \Omega ,\Lambda \right) $
is a spectral pair.
\end{corollary}

\begin{remark}
\label{Rem1.2}The simplest case of a disconnected $\Omega $ which has the
extension property, but which is not a spectral set, is in $d=1$, and we may
take $\Omega =\left\langle 0,1\right\rangle \cup \left\langle
2,4\right\rangle $, i.e., the union of two intervals with a doubling and
separation. The example was noted first in \cite{Fug74} and is based on the
simple observation that the polynomial $1+z^{2}+z^{3}$ has no roots $z$ on
the circle $\lvert z\rvert=1$.
\end{remark}

Some of the interest in spectral pairs derives from their connection to 
\emph{tilings.} A subset $\Omega \subset \mathbb{R}^{d}$ with nonzero
measure is said to \emph{tile }$\mathbb{R}^{d}$ if there is a set $L\subset %
\mathbb{R}^{d}$ such that the translates $\left\{ \Omega +l:l\in L\right\} $
cover $\mathbb{R}^{d}$ up to measure zero, and if the intersections 
\begin{equation}
\left( \Omega +l\right) \cap \left( \Omega +l^{\prime }\right) \text{\quad
for }l\ne l^{\prime }\text{ in }L  \label{eq10}
\end{equation}
have measure zero. We will call $\left( \Omega ,L\right) $ a \emph{tiling
pair} and we will say that $L$ is a \emph{set of translations}. The
\emph{Spectral-Set conjecture}
(see \cite{Fug74,Jor82,Ped87,
LaWa96c,LaWa97a}) states:

\begin{conjecture}
\label{Con1.3}Let $\Omega \subset \mathbb{R}^{d}$ have positive and finite
Lebesgue measure. Then $\Omega $ is a spectral set if and only if there is a
set $L$ of translations which make $\Omega $ tile $\mathbb{R}^{d}$.
\end{conjecture}

\begin{lemma}
\label{Lem1.4}If $\Omega =I^{d}$, then the zero-set for the function $%
F_{\Omega }$ in \textup{(\ref{eq4})} is 
\begin{equation}
\mathbf{Z}_{I^{d}}=\left\{ z\in \mathbb{C}^{d}\diagdown \left\{ 0\right\}
:\exists j\in \left\{ 1,\dots ,d\right\} \;\mathrm{s.t.}\;z_{j}\in \mathbb{Z}%
\diagdown \left\{ 0\right\} \right\} .  \label{eq11}
\end{equation}
\end{lemma}

\begin{proof}%
The function $F_{I^{d}}\left( \,\cdot \,\right) $ factors as follows. 
\begin{equation}
F_{I^{d}}\left( z\right) =\prod_{j=1}^{d}e^{i\pi z_{j}}\frac{\sin \pi z_{j}}{%
\pi z_{j}}  \label{eq12}
\end{equation}
for $z=\left( z_{1},\dots ,z_{d}\right) \in \mathbb{C}^{d}$, with the
interpretation that the function $z\mapsto \frac{\sin \pi z}{\pi z}$ is $1$
when $z=0$ in $\mathbb{C}$.%
\end{proof}%

\begin{remark}
\label{Rem1.7}What is special about $\mathbf{Z}_{\Omega }$ for $\Omega =I^{d}
$, as opposed to the general form of $\Omega $, is that $\mathbf{Z}%
_{I^{d}}\cup \left\{ 0\right\} $ is the \emph{Cayley graph} of the group $%
\Gamma =\mathbb{Z}^{d}$ with generators 
\begin{equation*}
S=\left\{ \left( \pm 1,0,\dots ,0\right) ,\dots ,\left( 0,\dots ,\pm
1,0,\dots ,0\right) ,\dots ,\left( 0,\dots ,0,\pm 1\right) \right\} .
\end{equation*}
We recall from \cite[Chapter 10]{BKS} the definition of the Cayley graph $%
G\left( \Gamma ,S\right) $ of a discrete group $\Gamma $ with generators $S$%
, $e\notin S$. When $\Gamma ,S$ are given, $G\left( \Gamma ,S\right) $ is
the graph with vertex set $\Gamma $ in which two vertices $\gamma
_{1},\gamma _{2}$ are the two ends of an edge iff $\gamma _{1}^{-1}\gamma
_{2}^{{}}\in S$. This gives a non-oriented graph, without any loop or
multiple edge.
\end{remark}

\section{Two dimensions\label{S2}}

We begin with the following simple observation in one dimension for $\Omega
=I=\left[ 0,1\right\rangle $. (For details, see \cite{JoPe92,ReSi}.)

\begin{proposition}
\label{Pro2.1}The only subsets $\Lambda \subset \mathbb{R}$ such that $%
\left( I,\Lambda \right) $ is a spectral pair are the translates 
\begin{equation}
\Lambda _{\alpha }:=\alpha +\mathbb{Z}=\left\{ \alpha +n:n\in \mathbb{Z}%
\right\}   \label{eq13}
\end{equation}
where $\alpha $ is some fixed real number.
\end{proposition}

In two dimensions, the corresponding result is more subtle, but the
possibilities may still be enumerated as follows:

\begin{theorem}
\label{Thm2.2}\textup{(\cite{JoPe99})} The only subsets $\Lambda \subset \mathbb{R}^{2}$ such that $%
\left( I^{2},\Lambda \right) $ is a spectral pair must belong to either one
or the other of the two classes, indexed by a number $\alpha $, and a
sequence $\left\{ \beta _{m}\in \left[ 0,1\right\rangle :m\in \mathbb{Z}%
\right\} $, where 
\begin{align}
\Lambda & =\left\{ 
\begin{pmatrix}
                   \alpha +m \\
                   \beta _{m}+n
                   \end{pmatrix}
                   %
:m,n\in \mathbb{Z}\right\}   \label{eq14} \\
\intertext{or}%
\Lambda & =\left\{ 
\begin{pmatrix}
                   \beta _{n}+m \\
                   \alpha +n
                   \end{pmatrix}
                   %
:m,n\in \mathbb{Z}\right\} .  \label{eq15}
\end{align}
Each of the two types occurs as the spectrum of a pair for the cube $I^{2}$,
and each of the sets $\Lambda $ as specified is a set of translation vectors
which produces a tiling of $\mathbb{R}^{2}$ by the cube $I^{2}$.
\end{theorem}

\begin{proof}%
See \cite{JoPe99} for details. The following
are some remarks of relevance to the general extension problem for operators.

The assertion in the theorem about $\Lambda $-translations tiling the plane
with $I^{2}$ is also clear from (\ref{eq14})--(\ref{eq15}), and it is
illustrated graphically in Figures \ref{tiling1} and \ref{tiling2}.

\TeXButton{figure1}{\begin{figure}[tbp]
\begin{center}
\includegraphics[height=4.0257in,width=3.3036in]{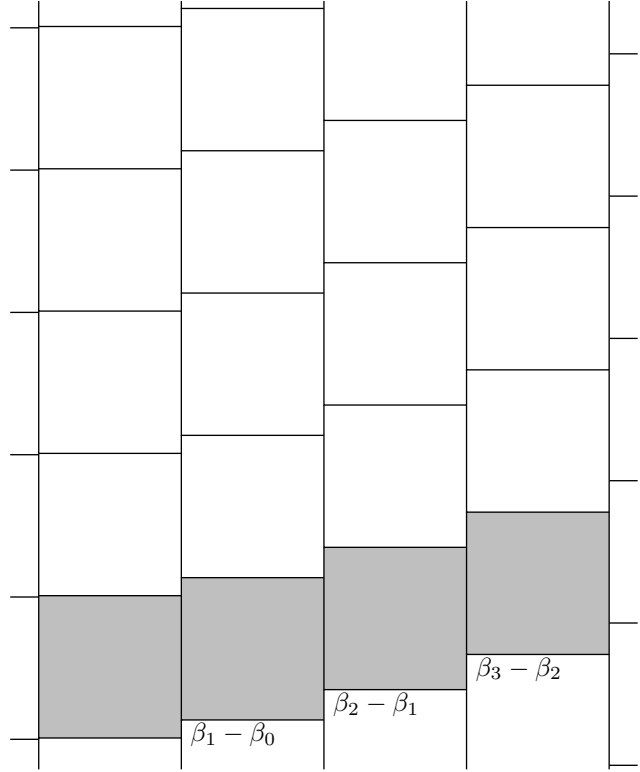}
\end{center}
\caption{\protect\figlabone Illustrating tiling with (\ref{eq14})}
\label{tiling1}
\end{figure}
}

\TeXButton{figure2}{\begin{figure}[tbp]
\begin{center}
\includegraphics[height=3.3036in,width=4.0257in]{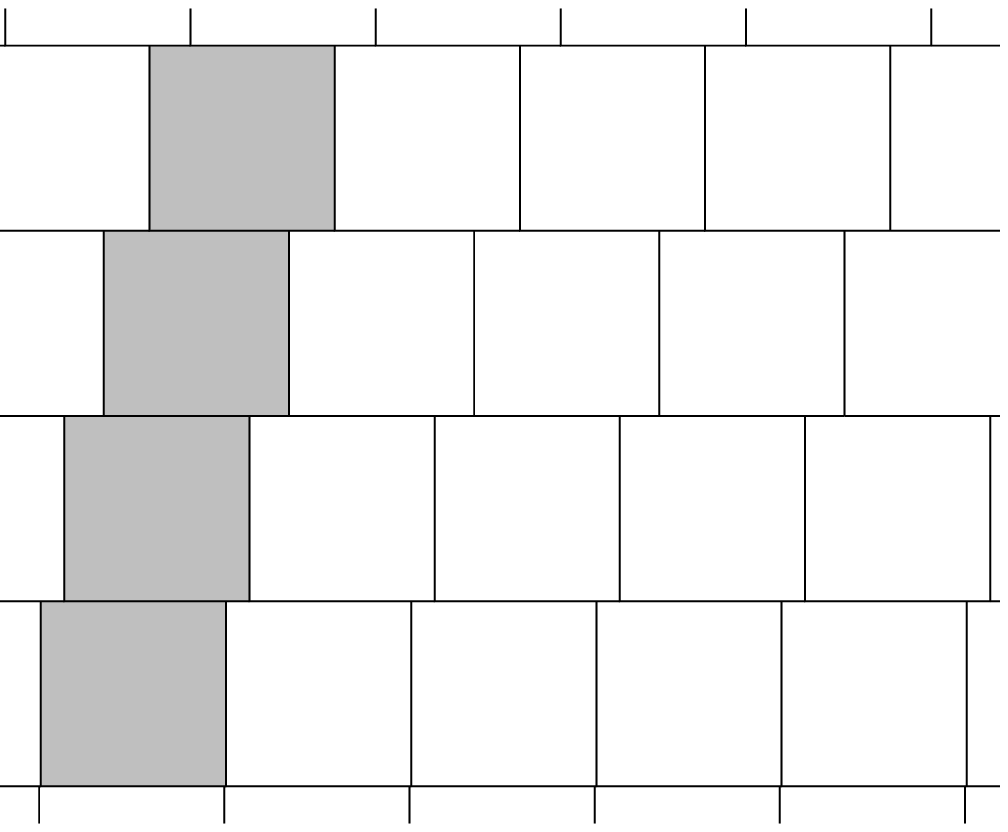}
\end{center}
\caption{\protect\figlabtwo Illustrating tiling with (\ref{eq15})}
\label{tiling2}
\end{figure}
}

It is clear that the pattern (\ref{eq14}) for $d=2$ continues to higher
dimensions as follows: 
\begin{equation}
\begin{pmatrix}
                   \alpha + k_{1} \\
                   \beta \left( k_{1}\right) +k_{2} \\
                   \gamma \left( k_{1},k_{2}\right) +k_{3} \\
                   \vdots  \\
                   \zeta \left( k_{1},k_{2},\dots ,k_{d-1}\right) +k_{d}
                   \end{pmatrix}
                   %
\label{eq17}
\end{equation}
with $k_{1},k_{2},\dots ,k_{d}\in \mathbb{Z}$, and 
\begin{align*}
\beta & \colon \mathbb{Z}\longrightarrow \left[ 0,1\right\rangle , \\
\gamma & \colon \mathbb{Z}\times \mathbb{Z}\longrightarrow \left[
0,1\right\rangle , \\
& \qquad \vdots \\
\zeta & \colon \mathbb{Z}^{d-1}\longrightarrow \left[ 0,1\right\rangle .
\end{align*}
Of course, then there are the obvious modifications of those cases resulting
from permutation of the $d$ coordinates; but the assertion is that, when $%
d\ge 10$, these configurations do \emph{not} suffice for cataloguing all the
possible spectra $\Lambda $ which turn $\left( I^{d},\Lambda \right) $ into
an $\mathbb{R}^{d}$-spectral pair.

We now turn to the non-trivial spectral-theoretic content of the conclusion
of the theorem. We claim that the two cases (\ref{eq14})--(\ref{eq15})
suffice when $d=2$. Note that the sequence $\beta \colon \mathbb{Z}%
\rightarrow \left[ 0,1\right\rangle $ is completely arbitrary.

We will show in Theorem \ref{Thm4.1} below that, up to a single translation
in the plane, the possibilities for the coordinates of points in a spectrum $%
\Lambda $ for $I^{2}$ are given by two sequences $\xi _{m}$, $\eta _{n}$
satisfying the following two cocycle relations: 
\begin{align}
\left( e^{i\xi _{m+k}}-e^{i\xi _{m}}\right) \left( 1-e^{i\eta _{n}}\right) &
=0  \label{eq18} \\
\intertext{and}%
\left( e^{i\eta _{n+l}}-e^{i\eta _{n}}\right) \left( 1-e^{i\xi _{m}}\right)
& =0  \label{eq19}
\end{align}
as identities in $m,n\in \mathbb{Z}$, and $k,l\in \mathbb{Z}\diagdown
\left\{ 0\right\} $. Note that the respective sequences are determined from
this only up to $2\pi \mathbb{Z}$ at each coordinate place.

Simple algebra shows that the two identities (\ref{eq18})--(\ref{eq19})
imply the following single identity 
\begin{equation}
\left( 1-e^{i\xi _{m+k}}\right) \left( 1-e^{i\eta _{n}}\right) =\left(
1-e^{i\xi _{m}}\right) \left( 1-e^{i\eta _{n+l}}\right)  \label{eq20}
\end{equation}
again for all $m,n\in \mathbb{Z}$ and $k,l\in \mathbb{Z}\diagdown \left\{
0\right\} $. But it follows from (\ref{eq20}) that at least one of the two
sequences, $1-e^{i\xi _{m}}$ or $1-e^{i\eta _{n}}$, must then vanish
identically. This yields the connection to the two cases for $\Lambda $
stated in (\ref{eq14})--(\ref{eq15}) of the theorem.

Hence the result giving two classes for $\Lambda $ in Theorem \ref{Thm2.2}
may be derived from our more general result in Section \ref{S4}.

The proof sketch of Theorem \ref{Thm2.2} is completed for now,
but details will be resumed in Section \ref{S4} below.%
\end{proof}%

\section{Operator extensions\label{S3}}

We saw in Theorem \ref{Thm1.1} that in some cases the existence problem for
spectral pairs, i.e., the question of when some given open subset $\Omega $
in $\mathbb{R}^{d}$ has an orthogonal basis $\left\{ e_{\lambda }:\lambda
\in \Lambda \right\} $ in $\mathcal{L}^{2}\left( \Omega \right) $ for some
set $\Lambda $ in $\mathbb{R}^{d}$, may be reformulated as a problem about
existence of commuting self-adjoint extensions of the operators $\left\{ 
\frac{1}{i}\frac{\partial \;}{\partial x_{j}}:j=1,\dots ,d\right\} $ with
common (dense) domain $C_{c}^{\infty }\left( \Omega \right) $ in $\mathcal{L}%
^{2}\left( \Omega \right) $. Suppose for the moment that $\Omega
=\left\langle 0,1\right\rangle \times \Omega _{2}$ where $\Omega _{2}$ is
some subset in $\mathbb{R}^{d-1}$ of finite positive $\left( d-1\right) $%
-dimensional Lebesgue measure. We then have the following classification of
the self-adjoint extensions $H$ of $\frac{1}{i}\frac{\partial \;}{\partial
x_{1}}$.

\begin{theorem}
\label{Lem3.1}The symmetric operator $\frac{1}{i}\frac{\partial \;}{\partial
x_{1}}$ in $\mathcal{L}^{2}\left( \left\langle 0,1\right\rangle \times
\Omega _{2}\right) $ with dense domain $\mathcal{D}$ consisting of $\varphi
\in \mathcal{L}^{2}\left( \left\langle 0,1\right\rangle \times \Omega
_{2}\right) $ such that $\varphi \left( \,\cdot \,,y\right) \in
C_{c}^{\infty }\left( \left\langle 0,1\right\rangle \right) $ for all $y\in
\Omega _{2}$, has self-adjoint extensions indexed by unitary operators $V$
in $\mathcal{L}^{2}\left( \Omega _{2}\right) $ in such a way that the 
\textup{(}unique\textup{)} extension $H_{V}$ is determined by its core
domain being of the form 
\begin{equation}
\mathcal{D}_{V}=\left\{ \varphi \left( x_{1},y\right) +e^{x_{1}}h\left(
y\right) +e^{1-x_{1}}\left( Vh\right) \left( y\right) :\varphi \in \mathcal{D%
},\,h\in \mathcal{L}^{2}\left( \Omega _{2}\right) \right\}   \label{eq22}
\end{equation}
and 
\begin{equation}
iH_{V}\left( \varphi +e^{x_{1}}h+e^{1-x_{1}}Vh\right) =\frac{\partial
\varphi }{\partial x_{1}}+e^{x_{1}}h-e^{1-x_{1}}Vh,  \label{eq23}
\end{equation}
for $\varphi \in \mathcal{D}$ and $h\in \mathcal{L}^{2}\left( \Omega
_{2}\right) $. We shall interpret the implicit boundary condition dictating
some extension $H_{V}$ as 
\begin{equation}
f\left( 1,\,\cdot \,\right) =U_{V}\left( f\left( 0,\,\cdot \,\right) \right)
,  \label{eq24}
\end{equation}
$f\in \mathcal{D}_{V}$ where the partial isometry $U_{V}$ is given by 
\begin{equation}
W_{V}=\left( eI+V\right) \left( I+eV\right) ^{-1},
\qquad U_{V}=\exp W_{V}.  \label{eq25}
\end{equation}
Conversely, $V$ may be calculated from $U_{V}$ by 
\begin{equation}
V=\left( I-eW_{V}\right) ^{-1}\left( W_{V}-eI\right) ,  \label{eq26}
\end{equation}
and in each case, the fractional linear transform, and its inverse, are well
defined.
\end{theorem}

\begin{proof}%
The proof is based on von Neumann's deficiency-space analysis of
self-adjoint extensions, and we refer to \cite{vNeu29}, \cite{ReSi}, and 
\cite{
Jor79} for background material on the theory of operator
extensions. If $S$ is a symmetric operator with dense domain $\mathcal{D}$
in a Hilbert space $\mathcal{H}$, then it has self-adjoint extensions if and
only if the two spaces 
\begin{equation}
\left( \left( iI\pm S\right) \mathcal{D}\right) ^{\perp }=:\mathcal{D}_{\pm }
\label{eq27}
\end{equation}
have the same dimensions. In that case, the corresponding extensions are
given by \emph{partial isometries} between the respective defect spaces $%
\mathcal{D}_{+}$ and $\mathcal{D}_{-}$ (see \cite{vNeu29}, \cite{ReSi}, or 
\cite{DS2}).
For convenience, we have chosen
a slightly different ``normalization'' in
our treatment of the Cayley
transform (\ref{eq25}) and its inverse (\ref{eq26}).
We did not normalize the
functions $e^{x_{1}}$ and $e^{1-x_{1}}$ in
the defect spaces. They have
$\mathcal{L}^{2}\left( I\right) $-norm equal to
$\left( \frac{e^{2}-1}{2}\right) ^{\frac{1}{2}}$. The fact that
$U_{V}$ in (\ref{eq25}) then defines a
partial isometry as claimed
amounts to the identities:
\[
\begin{minipage}[t]{\displayboxwidth}\raggedright
If $\psi \left( x_{1},y\right) =e^{x_{1}}h\left( y\right) 
+e^{1-x_{1}}\left( Vh\right) \left( y\right) $
as in (\ref{eq22}), then
\[
\psi \left( 1,y\right) =eh\left( y\right) 
+\left( Vh\right) \left( y\right) =\left( eI+V\right) h,
\]
and
\[
\psi \left( 0,y\right) =h\left( y\right) 
+eVh\left( y\right) =\left( I+eV\right) h.
\]
\end{minipage}
\]
This means that the vectors in the 
domain (\ref{eq22}) are given by the boundary
conditions (\ref{eq24}) which in turn determine
the unitary one-parameter group
\[
\mathcal{U}_{V}\left( t\right) :=\exp \left( itH_{V}\right) ,\qquad t\in \mathbb{R}.
\]
This group is defined from (\ref{eq24})
by using translation modulo $\mathbb{Z}$
in the $x_{1}$-variable. Then the
operator $U_{V}$ in (\ref{eq25}) is
used in defining the representation
$\mathbb{R}\ni t\mapsto \mathcal{U}_{V}\left( t\right) $ via
induction from $\mathbb{Z}$.

If $V\colon \mathcal{D}_{+}\rightarrow \mathcal{D}_{-}$ is a
partial isometry, then the domain of the corresponding extension $H$ ($%
H=H_{V}$) is 
\begin{equation*}
\left\{ \varphi +h_{+}+Vh_{+}:\varphi \in \mathcal{D},\,h_{+}\in \mathcal{D}%
_{+}\right\} 
\end{equation*}
and 
\begin{equation}
iH_{V}\left( \varphi +h_{+}+Vh_{+}\right) =iS\varphi +h_{+}-Vh_{+}.
\label{eq28}
\end{equation}
It follows that the lemma amounts to an identification of the \emph{defect
spaces} $\mathcal{D}_{\pm }$ when the symmetric operator is as specified.
When the variables in $\Omega =\left\langle 0,1\right\rangle \times \Omega
_{2}$ are separated as $\left( x_{1},y\right) $, $0<x_{1}<1$, $y=\left(
x_{2},\dots ,x_{d}\right) \in \Omega _{2}$, then vectors $h_{\pm }\in 
\mathcal{D}_{\pm }$ are precisely the solutions to 
\begin{equation}
S^{*}h_{\pm }=\pm ih_{\pm }.  \label{eq29}
\end{equation}
This amounts to solving 
\begin{equation*}
\frac{\partial \;}{\partial x_{1}}h_{\pm }\left( x_{1},y\right) =\pm h_{\pm
}\left( x_{1},y\right) 
\end{equation*}
in the sense of distributions, but with the restrictions $h_{\pm }\in 
\mathcal{L}^{2}\left( \left\langle 0,1\right\rangle \times \Omega
_{2}\right) $. The result of the lemma then follows from von Neumann's
characterization. If the minimal operator is not closed at the outset, then
the resulting self-adjoint extension comes from passing to the operator
closure in the formulas (\ref{eq23}) and (\ref{eq28}).%
\end{proof}%

\begin{corollary}
\label{Cor3.2}Let $V$ be a unitary operator in $\mathcal{L}^{2}\left( \Omega
_{2}\right) $ and let $H_{V}$ be the self-adjoint extension operator
described in Theorem \textup{\ref{Lem3.1}} in \textup{(\ref{eq23})--(\ref{eq24}%
).} Then $H_{V}$ generates a unitary one-parameter group $\left\{
U_{V}\left( t\right) :t\in \mathbb{R}\right\} $ in $\mathcal{L}^{2}\left(
\left\langle 0,1\right\rangle \times \Omega _{2}\right) $ which may be
realized \textup{(}up to unitary equivalence\textup{)} in the Hilbert space $%
\mathcal{H}_{V}$ of measurable functions $f\colon \mathbb{R}\rightarrow 
\mathcal{L}^{2}\left( \Omega _{2}\right) $, satisfying 
\begin{equation}
f\left( x_{1}+1\right) =U_{V}\left( f\left( x_{1}\right) \right) ,
\label{eq30}
\end{equation}
for all $x_{1}\in \mathbb{R}$, where $U_{V}$ is the operator from \textup{(%
\ref{eq25})} in Theorem \textup{\ref{Lem3.1}}, and the norm on
$\mathcal{H}_{V}$ is defined by 
\begin{equation*}
\left\| f\right\| _{\mathcal{H}_{V}}^{2}=\int_{0}^{1}\left\| f\left(
x_{1}\right) \right\| _{\mathcal{L}^{2}\left( \Omega _{2}\right)
}^{2}\,dx_{1}.
\end{equation*}
In this space the group $U_{V}\left( t\right) \colon \mathcal{H}%
_{V}\rightarrow \mathcal{H}_{V}$ is given by 
\begin{equation}
\left( U_{V}\left( t\right) f\right) \left( x_{1}\right) =f\left(
x_{1}-t\right) ,\text{\quad for }x_{1},t\in \mathbb{R}.  \label{eq31}
\end{equation}
The unitary isomorphism of $\mathcal{H}_{V}$ onto $\mathcal{L}^{2}\left(
\left\langle 0,1\right\rangle \times \Omega _{2}\right) =\mathcal{L}%
^{2}\left( \left\langle 0,1\right\rangle ,\mathcal{L}^{2}\left( \Omega
_{2}\right) \right) $ is simply the restriction to $\left\langle
0,1\right\rangle $ in the $x_{1}$-variable. Finally, if $U_{V}\left(
t\right) $ is computed in $\mathcal{L}^{2}\left( \left\langle
0,1\right\rangle \times \Omega _{2}\right) $, the formula is 
\begin{equation}
\left( U_{V}\left( t\right) f\right) f\left( x_{1},\,\cdot \,\right) =%
\begin{cases} f\left( x_{1}-t,\,\cdot \,\right) & \text{ if }0\le t<x_{1}<1,
\\ U_{V}\left( f\left( x_{1}-t,\,\cdot \,\right) \right) & \text{ if
}0<x_{1}\le t\le 1. \end{cases}  \label{eq32}
\end{equation}
\end{corollary}

\begin{proof}%
The realization on the space $\mathcal{H}_{V}$ is the interpretation of $%
U_{V}$ as a unitary representation of the group $\mathbb{R}$ which is
induced from the subgroup $\mathbb{Z}$ via formula (\ref{eq31}).
The advantage of this viewpoint is that
the spectral resolution of the unitary operator $U_{V}$ leads directly to an
associated direct integral decomposition for the unitary one-parameter group 
$\left\{ U_{V}\left( t\right) :t\in \mathbb{R}\right\} $ which is generated
by the extension operator $H_{V}$.%
\end{proof}%

When the corollary is applied to $\mathcal{L}^{2}\left( I\times I\right) $
from Section \ref{S2} we note that the respective unitary one-parameter
groups, $U_{x}\left( s\right) $ and $U_{y}\left( t\right) $, on $\mathcal{L}%
^{2}\left( I^{2}\right) $ which are generated by self-adjoint extension
operators of $\frac{1}{i}\frac{\partial \;}{\partial x}$ and $\frac{1}{i}%
\frac{\partial \;}{\partial y}$ with domain $C_{c}^{\infty }\left(
I^{2}\right) $, are induced representations in the sense of (\ref{eq30})--(%
\ref{eq31}). For the extensions of $\frac{1}{i}\frac{\partial \;}{\partial x}
$, the boundary-unitary from (\ref{eq30}) is acting on $\mathcal{L}%
^{2}\left( \left\{ 0<y<1\right\} \right) $. But we shall view it as a
unitary operator in $\mathcal{L}^{2}\left( I\times I\right) =\mathcal{L}%
^{2}\left( I_{x}\right) \otimes \mathcal{L}^{2}\left( I_{y}\right) $ via $%
U\leftrightarrow I\otimes U_{2}$ with $U_{2}$ acting in the $y$-variable. A
similar observation applies to the unitary one-parameter group $\left\{
U_{y}\left( t\right) :t\in \mathbb{R}\right\} $ acting on $\mathcal{L}%
^{2}\left( I^{2}\right) $ and generated by one of the self-adjoint
extensions of $\frac{1}{i}\frac{\partial \;}{\partial y}$. Hence the
boundary conditions for $\left\{ U_{x}\left( s\right) :s\in \mathbb{R}%
\right\} $ are given by a unitary $U\simeq I\otimes U_{2}$ with $U_{2}$
acting in the second variable, while those of $\left\{ U_{y}\left( t\right)
:t\in \mathbb{R}\right\} $ are determined by a second unitary operator $V$
in $\mathcal{L}^{2}\left( I^{2}\right) $, now of the form $V\leftrightarrow
V_{1}\otimes I$ with $V_{1}$ acting in the first variable of $\mathcal{L}%
^{2}\left( I\times I\right) $.

With this terminology we have the following preliminary result for the
square $I^{2}$ in the plane.

\begin{theorem}
\label{Thm3.3}Let $U_{x}\left( s\right) $ be the unitary one-parameter group
on $\mathcal{L}^{2}\left( I\times I\right) $, and let $U_{2}$ be the
corresponding unitary boundary operator acting in the second variable $y$.
Then $U_{2}$ commutes with the phase-periodic translation in the $y$%
-variable for a phase angle $\beta $ if and only if there is a real-valued
sequence $\left\{ \varphi _{n}:n\in \mathbb{Z}\right\} $ such that 
\begin{equation}
U_{x}\left( s\right) e_{m+\varphi _{n}}\otimes e_{n+\beta }=e^{i2\pi \left(
m+\varphi _{n}\right) s}e_{m+\varphi _{n}}\otimes e_{n+\beta }  \label{eq33}
\end{equation}
for all $s\in \mathbb{R}$ and $m,n\in \mathbb{Z}$, where for $\left( \xi
,\eta \right) \in \mathbb{R}^{2}$, $e_{\xi }\otimes e_{\eta }\left(
x,y\right) =e_{\xi }\left( x\right) e_{\eta }\left( y\right) =e^{i2\pi
\left( \xi x+\eta y\right) }$, restricted to $\left( x,y\right) \in I^{2}$.
\end{theorem}

\begin{proof}%
Recall that some fixed unitary one-parameter group $\left\{ U_{x}\left(
s\right) :s\in \mathbb{R}\right\} $ on $\mathcal{L}^{2}\left( I\times
I\right) $ is determined uniquely by the corresponding boundary operator $%
I\otimes U_{2}$. But it follows from Proposition \ref{Pro2.1} that $U_{2}$
satisfies the commutativity property of the theorem if and only if it is
diagonalized by the basis functions $\left\{ e_{n+\beta }:n\in \mathbb{Z}%
\right\} $ in $\mathcal{L}^{2}\left( I_{y}\right) $ for some $\beta \in
\left[ 0,1\right\rangle $, i.e., if, for some sequence $\varphi _{n}$, 
\begin{equation}
\left( U_{2}e_{n+\beta }\right) \left( y\right) =e^{i2\pi \varphi
_{n}}e_{n+\beta }\left( y\right) .  \label{eq3.13}
\end{equation}
But, according to Corollary \ref{Cor3.2}, this means that $U_{x}\left(
s\right) $ as an induced representation decomposes accordingly, which is to
say that the basis vectors $e_{m+\varphi _{n}}\otimes e_{n+\beta }$
simultaneously diagonalize each operator $U_{x}\left( s\right) $ as stated
in formula (\ref{eq33}).%
\end{proof}%

\begin{remark}
For more details on the operator-theoretic approach to spectrum and to
tiles, we refer to \cite{Jor87b,
Jor89b,
Ped87,Ped96}
\end{remark}

\section{Cocycles in two dimensions\label{S4}}

In this section, we continue with the self-adjoint extensions of the two
commuting minimal operators $\frac{1}{i}\frac{\partial \;}{\partial x}$ and $%
\frac{1}{i}\frac{\partial \;}{\partial y}$ with common dense domain $%
C_{c}^{\infty }\left( I^{2}\right) $ in $\mathcal{L}^{2}\left( I^{2}\right) $%
.

\begin{theorem}
\label{Thm4.1}Consider two commuting unitary one-parameter groups $%
U_{x}\left( s\right) $ and $U_{y}\left( t\right) $ on $\mathcal{L}^{2}\left(
I\times I\right) $ with respective boundary unitaries $U_{2}$ and $V_{1}$.
Then:

\begin{enumerate}
\item  \label{Thm4.1(1)}Either $U_{2}$ is of the form $aI_{\mathcal{L}%
^{2}\left( I_{y}\right) }$ for a scalar $a$, or else $V_{1}$ commutes with
periodic translation in the $x$-variable.

\item  \label{Thm4.1(2)}Either $V_{1}$ is of the form $bI_{\mathcal{L}%
^{2}\left( I_{x}\right) }$ for some scalar $b$, or else $U_{2}$ commutes
with periodic translation in the $y$-variable.

\item  \label{Thm4.1(3)}In case $U_{2}=e^{i2\pi \alpha }I_{\mathcal{L}%
^{2}\left( I_{y}\right) }$, then 
\begin{equation}
U_{x}\left( s\right) \left( e_{\alpha +m}\otimes g\right) =e^{i2\pi \left(
\alpha +m\right) s}e_{\alpha +m}\otimes g  \label{eq34}
\end{equation}
for all $m\in \mathbb{Z}$ and $g\in \mathcal{L}^{2}\left( I_{y}\right) $.

\item  \label{Thm4.1(4)}In case $V_{1}=e^{i2\pi \beta }I_{\mathcal{L}%
^{2}\left( I_{x}\right) }$, then 
\begin{equation}
U_{y}\left( t\right) \left( f\otimes e_{\beta +n}\right) =e^{i2\pi \left(
\beta +n\right) t}f\otimes e_{\beta +n}  \label{eq35}
\end{equation}
for all $f\in \mathcal{L}^{2}\left( I_{x}\right) $ and $n\in \mathbb{Z}$.
\end{enumerate}
\end{theorem}

\begin{remark}
\label{Rem4.2}It follows that the conclusion in Theorem \ref{Thm3.3} is
satisfied when the two one-parameter groups commute, i.e., when 
\begin{equation}
U_{x}\left( s\right) U_{y}\left( t\right) =U_{y}\left( t\right) U_{x}\left(
s\right)   \label{eq36}
\end{equation}
is assumed, $s,t\in \mathbb{R}$. Specifically, it will then always be the
case that $U_{2}$ commutes with some phase-periodic translation in the $y$%
-variable, while $V_{1}$ commutes with some (possibly different)
phase-periodic translation in the $x$-variable. (Also note that (\ref{eq36})
is a reformulation of (\ref{eq9}) in the case $d=2$. Furthermore (\ref{eq36}%
) signifies the presence of a unitary representation of $\mathbb{R}^{2}.$)
\end{remark}

\begin{proof}[Proof of Theorem \textup{\ref{Thm4.1}}]%
When the two one-parameter groups $U_{x}\left( s\right) $ and $U_{y}\left(
t\right) $ are written in the form (\ref{eq32}) from Corollary \ref{Cor3.2},
then the alternatives in (\ref{eq32}) may be expanded as follows. Let $\tau
_{s}$ denote periodic translation in $\mathcal{L}^{2}\left( \left\langle
0,1\right\rangle \right) $, and let $P_{s}$ denote the projection of $%
\mathcal{L}^{2}\left( \left\langle 0,1\right\rangle \right) $ onto $\mathcal{%
L}^{2}\left( \left\langle 0,s\right\rangle \right) $, with $P_{s}^{\perp
}=I-P_{s}$ denoting then the projection onto the complement $\mathcal{L}%
^{2}\left( \left\langle s,1\right\rangle \right) $, for $s\in \left[
0,1\right] $. We have $P_{0}=0$ and $P_{1}=I_{\mathcal{L}^{2}\left(
\left\langle 0,1\right\rangle \right) }$. Then from (\ref{eq32}) we get 
\begin{align}
U_{x}\left( s\right) & =\tau _{s}^{{}}P_{s}^{\perp }\otimes I+\tau
_{s}P_{s}\otimes U_{2}  \label{eq37} \\
\intertext{and}%
U_{y}\left( t\right) & =I\otimes \tau _{t}^{{}}P_{t}^{\perp }+V_{1}\otimes
\tau _{t}P_{t}.  \label{eq38}
\end{align}
The assumed commutativity (\ref{eq36}) then takes the form: 
\begin{multline*}
\tau _{s}^{{}}P_{s}^{\perp }V_{1}^{{}}\otimes \tau _{t}P_{t}+\tau
_{s}P_{s}\otimes U_{2}^{{}}\tau _{t}^{{}}P_{t}^{\perp }+\tau
_{s}P_{s}V_{1}\otimes U_{2}\tau _{t}P_{t} \\
=\tau _{s}P_{s}\otimes \tau _{t}^{{}}P_{t}^{\perp }U_{2}^{{}}+V_{1}^{{}}\tau
_{s}^{{}}P_{s}^{\perp }\otimes \tau _{t}P_{t}+V_{1}\tau _{s}P_{s}\otimes
\tau _{t}P_{t}U_{2}.
\end{multline*}
If $V_{1}$ is not a scalar times $I_{\mathcal{L}^{2}\left( I_{x}\right) }$
then two terms on either side are independent when evaluated on $f\otimes g$%
. Hence both $U_{2}^{{}}\tau _{t}^{{}}P_{t}^{\perp }=\tau
_{t}^{{}}P_{t}^{\perp }U_{2}^{{}}$ and $U_{2}\tau _{t}P_{t}=\tau
_{t}P_{t}U_{2}$ hold. Addition of these two identities yields $U_{2}\tau
_{t}=\tau _{t}U_{2}$ which is the commutativity of $U_{2}$ with periodic
translation.

If on the other hand $V_{1}$ is a scalar, then it follows from the argument
in Section \ref{S3} that (\ref{Thm4.1(4)}) must hold.

The two possibilities for the other boundary operator $U_{2}$ lead to cases (%
\ref{Thm4.1(1)}) and (\ref{Thm4.1(3)}) by symmetry.%
\end{proof}%

\begin{corollary}
\label{Cor4.3}Consider unitary one-parameter groups $U_{x}\left( s\right) $
and $U_{y}\left( t\right) $ as in Theorem \textup{\ref{Thm4.1}} and suppose
the corresponding boundary operators $U_{2}$ and $V_{1}$ diagonalize as
follows \textup{(}identities in $n,m\in \mathbb{Z}$\textup{):} 
\begin{align}
U_{2}e_{n+\beta }& =e^{i2\pi \alpha _{n}}e_{n+\beta }  \label{eq39} \\
\intertext{and}%
V_{1}e_{m+\alpha }& =e^{i2\pi \beta _{m}}e_{m+\alpha }  \label{eq40}
\end{align}
for some $\alpha ,\beta \in \mathbb{R}$. The sequences $\alpha _{n},\beta
_{m}$ will be chosen taking values in $\left[ 0,1\right\rangle $. Then the
commutativity \textup{(\ref{eq36})} for the two groups holds if and only if
the two sequences satisfy a certain cocycle property: Let $a_{n}:=e^{i2\pi
\alpha _{n}}$ and $b_{m}:=e^{i2\pi \beta _{m}}$. Then the two identities 
\begin{align}
\left( b_{m}-b_{m+k}\right) \left( 1-a_{n}\right) & =0,\quad m,n\in %
\mathbb{Z},\;k\in \mathbb{Z}\diagdown \left\{ 0\right\}   \label{eq41} \\
\intertext{and}%
\left( a_{n}-a_{n+l}\right) \left( 1-b_{m}\right) & =0,\quad m,n\in %
\mathbb{Z},\;l\in \mathbb{Z}\diagdown \left\{ 0\right\}   \label{eq42}
\end{align}
are equivalent to the commutativity \textup{(\ref{eq36})}. If commutativity
holds, we must have $\left( 1-a_{n}\right) \left( 1-b_{m}\right) \equiv 0$, $%
n,m\in \mathbb{Z}$. Hence we get a spectral pair with spectrum $\Lambda $
having one of the two forms 
\begin{align}
\left\{ 
\begin{pmatrix}
                   \alpha +m \\
                   n+\beta _{m}
                   \end{pmatrix}
                   %
:m,n\in \mathbb{Z}\right\} \text{\quad if }\alpha _{n}& \equiv 0,  \tag{i}
\label{Cor4.3(1)} \\
\intertext{or}%
\left\{ 
\begin{pmatrix}
                   m+\alpha _{n} \\
                   \beta +n
                   \end{pmatrix}
                   %
:m,n\in \mathbb{Z}\right\} \text{\quad if }\beta _{m}& \equiv 0.  \tag{ii}
\label{Cor4.3(2)}
\end{align}
\end{corollary}

The derivation of the two cocycle identities (\ref{eq41})--(\ref{eq42}) from
commutativity (\ref{eq36}) at the end of the proof is based on the following
corollary of independent interest:

\begin{corollary}
\label{Lem4.4}Let $U=I\otimes U_{2}$ and $V=V_{1}\otimes I$ be the
respective boundary operators of the one-parameter unitary groups $%
U_{x}\left( s\right) $ and $U_{y}\left( t\right) $ acting on $\mathcal{L}%
^{2}\left( I\times I\right) $. Then, if \textup{(\ref{eq39})--(\ref{eq40})}
hold for some $\alpha ,\beta $ and some sequences as specified, it follows
that the respective one-parameter groups may be expanded in the common basis 
$E\left( m,n\right) =E_{\left( \alpha ,\beta \right) }\left( m,n\right)
:=e_{m+\alpha }^{\left( 1\right) }\otimes e_{n+\beta }^{\left( 2\right) }$, $%
\left( m,n\right) \in \mathbb{Z}^{2}$, as follows: There are complex
sequences $\left\{ s_{k}\right\} _{k\in \mathbb{Z}}$ and $\left\{
t_{l}\right\} _{l\in \mathbb{Z}}$ so that, if we define 
\begin{gather}
s_{0}^{\perp }:=1-s_{0}^{{}},\quad t_{0}^{\perp }:=1-t_{0}^{{}}  \label{eq43} \\
\intertext{and}%
s_{k}^{\perp }:=-s_{k}^{{}}\text{\quad \textup{(}for }k\ne 0\text{\textup{%
),\quad }}t_{l}^{\perp }:=-t_{l}^{{}}\text{\quad \textup{(}for }l\ne 0\text{%
\textup{),}}  \label{eq43bis}
\end{gather}
then
\begin{align}
U_{x}\left( s\right) E\left( m,n\right) & =\sum_{k\in \mathbb{Z}}e^{i2\pi
\left( m+\alpha +k\right) s}\left( s_{k}^{\perp
}+s_{k}^{{}}a_{n}^{{}}\right) E\left( m+k,n\right)   \label{eq44} \\
\intertext{and}%
U_{y}\left( t\right) E\left( m,n\right) & =\sum_{l\in \mathbb{Z}}e^{i2\pi
\left( n+\beta +l\right) t}\left( t_{l}^{\perp }+t_{l}^{{}}b_{m}^{{}}\right)
E\left( m,n+l\right) .  \label{eq45}
\end{align}
The two one-parameter groups $U_{x}\left( s\right) $ and $U_{y}\left(
t\right) $ commute if and only if the cocycle identities \textup{(\ref{eq41}%
)--(\ref{eq42})} hold.
\end{corollary}

\begin{proof}%
Recall from (\ref{eq37})--(\ref{eq38}) that the two one-parameter groups are
expressed in terms of multiplication operators on $\mathcal{L}^{2}\left(
\left\langle 0,1\right\rangle \right) $ with the respective indicator
functions $\chi _{\left\langle 0,s\right\rangle }$ and $\chi _{\left\langle
0,t\right\rangle }$. The sequences (\ref{eq43})--(\ref{eq43bis}) are the
Fourier coefficients of these indicator functions, acting by multiplication
in $\mathcal{L}^{2}\left( I\right) $, and the relations (\ref{eq43})--(\ref
{eq43bis}) simply reflect the following two obvious identities, 
\begin{align*}
\chi _{\left\langle 0,s\right\rangle }+\chi _{\left[ s,1\right\rangle }& =1
\\
\intertext{and}%
\chi _{\left\langle 0,t\right\rangle }+\chi _{\left[ t,1\right\rangle }& =1,
\end{align*}
as functions on the unit interval. When the resulting formulas (\ref{eq44}%
)--(\ref{eq45}) are substituted into 
\begin{equation}
U_{x}\left( s\right) U_{y}\left( t\right) E\left( m,n\right) =U_{y}\left(
t\right) U_{x}\left( s\right) E\left( m,n\right)  \label{eq46}
\end{equation}
the equivalence to (\ref{eq41})--(\ref{eq42}) results.%
\end{proof}%

\section{\label{SNew5}Quasicrystals}

For the spectral pairs $\left( I^{d},\Lambda \right) $ in dimensions $d=2,3$%
, we noted that each of the candidates for spectrum $\Lambda $ tiles $%
\mathbb{R}^{d}$ with $\Lambda $-translates of $I^{d}$. (See Theorems \ref
{Thm2.2} and \ref{Thm3.3}.)
But reviewing formulas (\ref{eq14})--(\ref{eq15}) and (\ref{eq17}%
), and (\ref{eq50}) in the next section, for the possible sets 
$\Lambda $ which serve as $I^{d}$-spectrum, we find functions $\alpha ,\beta
,\dots $ on $\mathbb{Z}$ or $\mathbb{Z}^{k}$ which describe the particular
set $\Lambda $. Since all the candidates for $\Lambda $ make tilings, there
is a direct \emph{geometric} interpretation for these functions; but we note
in the present section that there is also a \emph{spectral-theoretic}
significance which derives from diffraction considerations of quasicrystals;
see \cite{Sen95}, \cite{Hof95}, and \cite{BoTa87}.

In this setting, diffractions show up as discrete components of the spectral
distribution 
\begin{equation*}
D_{\Lambda }\left( x\right) =\sum_{\lambda \in \Lambda }e_{\lambda }\left(
x\right) =\sum_{\lambda \in \Lambda }e^{i2\pi \lambda \cdot x}.
\end{equation*}
We say that a spectrum $\Lambda $ ($\subset \mathbb{R}^{d}$) has a \emph{%
diffraction pattern} if there is a pair $\left( M,c\right) $ where $M$ is a
subset of $\mathbb{R}^{d}$ and $c$ is a function (measuring intensity)
defined on $M$ such that 
\begin{equation*}
D_{\Lambda }\left( x\right) =\sum_{m\in M}c\left( \mu \right) \delta \left(
x-\mu \right) ,
\end{equation*}
i.e., the spectral distribution is a weighted sum of point-masses, supported
on some (discrete) subset $M$ in $\mathbb{R}^{d}$. Note that the
interpretation in both of the summations involving $D_{\Lambda }\left(
\,\cdot \,\right) $ is to be understood as Schwartz distributions; that is
if the respective sums are evaluated on a testing function $\varphi \in
C_{c}^{\infty }\left( \mathbb{R}^{d}\right) $, then the first sum yields $%
\sum_{\lambda \in \Lambda }\tilde{\varphi}\left( \lambda \right) $ where $%
\tilde{\varphi}\left( \lambda \right) =\int_{\mathbb{R}^{d}}e_{\lambda
}\left( x\right) \varphi \left( x\right) \,dx$, while the second sum is $%
\sum_{\mu \in M}c\left( \mu \right) \varphi \left( \mu \right) $. We also
note that, by the Poisson summation formula, the condition is satisfied if $%
\Lambda =\mathbb{Z}^{d}=M$, and the density (intensity) function $c$ is $%
c\equiv 1$ on $M$.

We shall also need the following definition: A function $\xi $ on $\mathbb{R}
$ is said to be \emph{quasi-periodic} if there are positive numbers $\omega
_{1},\dots ,\omega _{r}$, which are independent over $\mathbb{Q}$, and
functions $\xi _{1},\dots ,\xi _{r}$ such that $\xi _{j}$ has $\omega _{j}$
as period, and $\xi =\sum_{j=1}^{r}\xi _{j}$. The condition on $\xi _{j}$
amounts to the generalized Fourier expansion 
\begin{equation*}
\xi _{j}\left( x\right) =\sum_{n\in \mathbb{Z}}c_{j}\left( n\right) e^{i2\pi 
\frac{nx}{\omega _{j}}}.
\end{equation*}

In the following result we show that, if the functions which define a
spectrum $\Lambda $ for some $I^{d}$ are quasi-periodic, then it follows
that $\Lambda $ has a diffraction pattern. We will not state the result in
the widest generality as it will be clear that the idea in the simplest case
carries over to the variations in higher dimensions. Even for $d=2$, Theorem 
\ref{Thm2.2} shows that there are two classes of $\Lambda $ corresponding to
(\ref{eq14}) and (\ref{eq15}) respectively. In the following we will treat
only (\ref{eq14}), but the result applies to (\ref{eq15}) \emph{mutatis
mutandis.}

\begin{theorem}
\label{ThmNew5.1}Let 
\begin{equation*}
\Lambda =\left\{ 
\begin{pmatrix}
                   m \\
                   \beta \left( m\right) +n
                   \end{pmatrix}
                   %
:m,n\in \mathbb{Z}\right\} 
\end{equation*}
for some function $\beta \colon \mathbb{Z}\rightarrow \mathbb{R}$
and suppose $\beta $ extends to a function on $\mathbb{R}$
which is quasi-periodic with periods $\omega _{1},\dots ,\omega _{r}$,
independent over $\mathbb{Q}$. Then it follows that $\left( I^{2},\Lambda
\right) $ is a spectral pair with diffraction pattern; specifically, there
is a density function $c\colon \mathbb{Z}^{r}\times \mathbb{Z}\rightarrow %
\mathbb{C}$ such that 
\begin{equation*}
D_{\Lambda }\left( x,y\right) =\sum_{k\in \mathbb{Z}^{r}}\sum_{n\in %
\mathbb{Z}}c\left( k,n\right) \sum_{m\in \mathbb{Z}}\delta \left( x-%
\smash{\sum_{i=1}^{r}}%
\frac{k_{i}}{\omega _{i}}-m\right) \otimes \delta \left( y-n\right) 
\end{equation*}
with the density $c\left( k_{1},\dots ,k_{r},n\right) $ derived from the
Bohr almost periodic Fourier expansion applied to $\beta $.
\end{theorem}

\begin{proof}%
Consider the formula $D_{\Lambda }\left( x,y\right)
=\sum_{m}\sum_{n}e^{i2\pi \left( mx+\left( \beta \left( m\right) +n\right)
y\right) }$ and expand the inside function, $m\mapsto e^{i2\pi \beta \left(
m\right) y}$ according to the quasi-periodicity assumption on $\beta $:
specifically, 
\begin{align*}
e^{i2\pi \beta \left( m\right) y}& =\prod_{j=1}^{r}e^{i2\pi \xi _{j}\left(
m\right) y} \\
& =\prod_{j=1}^{r}\sum_{k_{j}\in \mathbb{Z}}c^{\left( j\right) }\left(
k_{j}\right) e^{i2\pi \frac{mk_{j}}{\omega _{j}}} \\
& =\sum_{k_{1}\in \mathbb{Z}}\cdots \sum_{k_{r}\in \mathbb{Z}}c^{\left(
1\right) }\left( k_{1}\right) \cdots c^{\left( r\right) }\left( k_{r}\right)
e^{i2\pi m\sum_{j=1}^{r}\frac{k_{j}}{\omega _{j}}}\,.
\end{align*}
Setting $c\left( k\right) :=\prod_{j=1}^{r}c^{\left( j\right) }\left(
k_{j}\right) $ and using 
\begin{equation*}
\sum_{m\in \mathbb{Z}}e^{i2\pi m\left( x+\sum_{j=1}^{r}\frac{k_{j}}{\omega
_{j}}\right) }=\sum_{m\in \mathbb{Z}}\delta \left( x-%
\smash{\sum_{j=1}^{r}}%
\frac{k_{j}}{\omega _{j}}-m\right)
\end{equation*}
together with Poisson summation (also in the second variable) we arrive at
the desired formula.%
\end{proof}%

\section{\label{S5}Higher dimensions}

The following definitions help summarize the results for $d=2$: We say that
the one-parameter unitary groups on $\mathcal{L}^{2}\left( I\times I\right) $
generated by self-adjoint extensions of the respective partial derivatives $%
\frac{1}{i}\frac{\partial \;}{\partial x}$ and $\frac{1}{i}\frac{\partial \;%
}{\partial y}$ on $C_{c}^{\infty }\left( I\times I\right) $ are \emph{%
quasi-commuting} if the conditions (\ref{eq39})--(\ref{eq40}) hold. Recall
this means that the respective boundary operators commute with some
phase-periodic translation in the opposite variable. We then showed in
Theorem \ref{Thm4.1} that the commutativity property (\ref{eq36}), for the
unitary groups $U_{x}\left( s\right) $ and $U_{y}\left( t\right) $, implies 
\emph{quasi-commutativity.} Finally we showed in Corollary \ref{Cor4.3}
that, among the quasi-commuting extensions, those that in fact commute (in
the sense of (\ref{eq36})) are characterized by the two cocycle identities (%
\ref{eq41})--(\ref{eq42}).

It is clear that \emph{quasi-commutativity} can be defined analogously for $%
d>2$. It follows from Theorem \ref{Thm1.1} that commutativity of $d$
self-adjoint extensions of the respective partial derivatives $\left\{ \frac{%
1}{i}\frac{\partial \;}{\partial x_{j}}:j=1,\dots ,d\right\} $, on $%
C_{c}^{\infty }\left( I^{d}\right) \subset \mathcal{L}^{2}\left(
I^{d}\right) $, is equivalent to the \emph{spectral-pair} condition for $%
\left( I^{d},\Lambda \right) $. Moreover, if commuting self-adjoint
extensions exist (i.e., $\frac{1}{i}\frac{\partial \;}{\partial x_{j}}%
\subset H_{j}$, $H_{j}^{*}=H_{j}^{{}}$, $j=1,\dots ,d$), then we may take $%
\Lambda $ to be the joint spectrum of the family $\left\{ H_{j}\right\}
_{j=1}^{d}$. Conversely, commuting operators $H_{j}$ may easily be
associated with some spectrum $\Lambda $ in a spectral pair $\left(
I^{d},\Lambda \right) $. Hence, for $d=2$, our results in Section \ref{S4}
provide a complete classification of the commuting (and also the
quasi-commuting) self-adjoint extensions of $\left\{ \frac{1}{i}\frac{%
\partial \;}{\partial x_{j}}\right\} _{j=1}^{d}$.

In higher dimensions, we still have boundary operators corresponding to each
self-adjoint extension of the partials $\frac{1}{i}\frac{\partial \;}{%
\partial x_{j}}$ (on $C_{c}^{\infty }\left( I^{d}\right) \subset \mathcal{L}%
^{2}\left( I^{d}\right) $, $j=1,\dots ,d$), by Corollary \ref{Cor3.2}. If
for each $j$, $U_{j}\left( t\right) $ denotes the unitary one-parameter
group on $\mathcal{L}^{2}\left( I^{d}\right) $ generated by some
self-adjoint extension $H_{j}$, then Corollary \ref{Cor3.2} states that $%
U_{j}\left( t\right) $ is induced by some unitary operator $V_{j}$ acting in
the remaining variables $\left( x_{1},\dots ,x_{j-1},x_{j+1},\dots
,x_{d}\right) $ (i.e., with omission of the variable on the $j$'th place):
specifically, $U_{j}\left( t\right) =\limfunc{ind}_{\mathbb{Z}}^{\mathbb{R}%
}\left( V_{j}\right) $ as a representation of $\left( \mathbb{R},+\right) $;
or equivalently the domain of $H_{j}$ is, for each $j$, given by the
boundary condition 
\begin{equation*}
f\left( x_{1}\dots ,x_{j-1},1,x_{j+1},\dots ,x_{d}\right) =V_{j}\left(
f\left( x_{1}\dots ,x_{j-1},0,x_{j+1},\dots ,x_{d}\right) \right) .
\end{equation*}
(Note that the more precise interpretation of this set of boundary
conditions is given in formula (\ref{eq30}) of Corollary \ref{Cor3.2}. This
is the interpretation of the unitary one-parameter groups in the respective
coordinate directions as \emph{induced unitary representations} (see \cite
{Mac53,Mac62}), with induction $\mathbb{Z}\rightarrow \mathbb{R}$ for each
direction.) We say that a family of self-adjoint extension operators $H_{j}$%
, with corresponding boundary unitaries $V_{j}$, is \emph{quasi-commuting}
if there are phase angles $\alpha _{j}\in \left[ 0,1\right\rangle $, $%
j=1,\dots ,d$, such that each $V_{j}$ is diagonalized by 
\begin{equation}
e_{\alpha _{1}+n_{1}}^{\left( 1\right) }\otimes \dots \otimes e_{\alpha
_{j-1}+n_{j-1}}^{\left( j-1\right) }\otimes e_{\alpha
_{j+1}+n_{j+1}}^{\left( j+1\right) }\otimes \dots \otimes e_{\alpha
_{d}+n_{d}}^{\left( d\right) }  \label{eq47}
\end{equation}
as $\left( n_{1},\dots ,n_{d-1},n_{d+1},\dots ,n_{d}\right) $ vary over $%
\mathbb{Z}^{d-1}$; i.e., the lattice resulting from $\mathbb{Z}^{d}$ with
the $j$'th coordinate place omitted. It follows that the quasi-commutative
case is characterized by the phase angles $\alpha _{1},\dots ,\alpha _{d}$,
and by functions $v_{j}\colon \mathbb{Z}^{d-1}\rightarrow \mathbb{T}$ such
that, for $n=\left( n_{1},\dots ,%
\smash{\hat{\jmath }}%
,\dots ,n_{d}\right) $, $v_{j}\left( n\right) =v_{j}\left( n_{1},\dots ,%
\smash{\hat{\jmath }}%
,\dots ,n_{d}\right) $ is the eigenvalue of $V_{j}$ coresponding to the
eigenvector in (\ref{eq47}). (The notation $\left( n_{1},\dots ,%
\smash{\hat{\jmath }}%
,\dots ,n_{d}\right) $ means that the $j$'th place is omitted.)

\begin{theorem}
\label{Thm5.1}Let $\left\{ H_{j}\right\} _{j=1}^{d}$ be a family of
self-adjoint extensions of the respective partials $\frac{1}{i}\frac{%
\partial \;}{\partial x_{j}}$ \textup{(}$j=1,\dots ,d$, on $C_{c}^{\infty
}\left( I^{d}\right) \subset \mathcal{L}^{2}\left( I^{d}\right) $\textup{)},
which is assumed quasi-commutative with eigenvalue functions $v_{j}\left(
n_{1},\dots ,%
\smash{\hat{\jmath }}%
,\dots ,n_{d}\right) $ from $\mathbb{Z}^{d-1}$ to $\mathbb{T}$. Then the
extensions are commutative if and only if the following pair of cocycle
conditions is satisfied for all $j,k$ such that $1\le j<k\le d$, all $\left(
n_{1},\dots ,%
\smash{\hat{\jmath }}%
,\dots ,n_{d}\right) $, and all $l,m\in \mathbb{Z}\diagdown \left\{
0\right\} $: 
\begin{multline}
\left( v_{j}\left( n_{1},\dots ,%
\smash{\hat{\jmath }}%
,\dots ,n_{k}+l,\dots ,n_{d}\right) -v_{j}\left( n_{1},\dots ,%
\smash{\hat{\jmath }}%
,\dots ,n_{d}\right) \right)  \\
\times \left( 1-v_{k}\left( n_{1},\dots ,%
\smash{\hat{k}}%
,\dots ,n_{d}\right) \right) =0  \label{eq48}
\end{multline}
and 
\begin{multline}
\left( v_{k}\left( n_{1},\dots ,n_{j}+m,\dots ,%
\smash{\hat{k}}%
,\dots ,n_{d}\right) -v_{k}\left( n_{1},\dots ,%
\smash{\hat{k}}%
,\dots ,n_{d}\right) \right)  \\
\times \left( 1-v_{j}\left( n_{1},\dots ,%
\smash{\hat{\jmath }}%
,\dots ,n_{d}\right) \right) =0.  \label{eq49}
\end{multline}
\end{theorem}

\begin{proof}%
Since the commutativity for the one-parameter groups of unitary operators $%
U_{j}\left( t_{j}\right) $ may be stated for pairs, i.e., $U_{j}\left(
t_{j}\right) U_{k}\left( t_{k}\right) =U_{k}\left( t_{k}\right) U_{j}\left(
t_{j}\right) $, $j<k$, $t_{j}\in \mathbb{R}$, $t_{k}\in \mathbb{R}$, the
argument for the general case $d>2$ is the same as for $d=3$. To see this,
just use the formulas for the respective one-parameter groups which are
analogues to (\ref{eq44})--(\ref{eq45}) in the proof of Corollary \ref{Lem4.4}.
For $d=3$, we may introduce the leg-notation: $v_{1}\rightarrow v_{23}$, $%
v_{2}\rightarrow v_{13}$, $v_{3}\rightarrow v_{12}$. When evaluated at a
general point in $\mathbb{Z}^{3}$ of the form $\left( k,l,m\right) $, the
respective eigenvalues are: 
\begin{align}
v_{23}\left( l,m\right) & \text{\quad for }V_{23},  \tag{i}
\label{Thm5.1(a)(1)} \\
v_{13}\left( k,m\right) & \text{\quad for }V_{13},  \tag{ii}
\label{Thm5.1(a)(2)} \\
v_{12}\left( k,l\right) & \text{\quad for }V_{12}.  \tag{iii}
\label{Thm5.1(a)(3)}
\end{align}
Specifically, 
\begin{align}
V_{23}e_{\beta +l}^{\left( 2\right) }\otimes e_{\gamma +m}^{\left( 3\right)
}& =v_{23}\left( l,m\right) e_{\beta +l}^{\left( 2\right) }\otimes e_{\gamma
+m}^{\left( 3\right) },  \tag{i}  \label{Thm5.1(b)(1)} \\
V_{13}e_{\alpha +k}^{\left( 1\right) }\otimes e_{\gamma +m}^{\left( 3\right)
}& =v_{13}\left( k,m\right) e_{\alpha +k}^{\left( 1\right) }\otimes
e_{\gamma +m}^{\left( 3\right) },  \tag{ii}  \label{Thm5.1(b)(2)} \\
V_{12}e_{\alpha +k}^{\left( 1\right) }\otimes e_{\beta +l}^{\left( 2\right)
}& =v_{12}\left( k,l\right) e_{\alpha +k}^{\left( 1\right) }\otimes e_{\beta
+l}^{\left( 2\right) },  \tag{iii}  \label{Thm5.1(b)(3)}
\end{align}
where $\alpha ,\beta ,\gamma $ are the fixed phase angles from the
quasi-commutativity. Then the three pairs of cocycle identities from the
theorem are as follows: (\ref{Thm5.1(c)(1a)})--(\ref{Thm5.1(c)(1b)}), (\ref
{Thm5.1(c)(2a)})--(\ref{Thm5.1(c)(2b)}), and (\ref{Thm5.1(c)(3a)})--(\ref
{Thm5.1(c)(3b)}) below. The argument for the equivalence of commutativity
and the cocycle identities is essentially the same as the one used in the
proof of Corollary \ref{Lem4.4} above. The cocycle identities for $d=3$ are: 
\begin{align}
\left( v_{13}\left( k,m\right) -v_{13}\left( k+n_{1},m\right) \right) \left(
1-v_{23}\left( l,m\right) \right) & =0  \tag{i\thinspace a}
\label{Thm5.1(c)(1a)} \\
\left( v_{23}\left( l,m\right) -v_{23}\left( l+n_{2},m\right) \right) \left(
1-v_{13}\left( k,m\right) \right) & =0,  \tag{i\thinspace b}
\label{Thm5.1(c)(1b)}
\end{align}
\begin{align}
\left( v_{12}\left( k,l\right) -v_{12}\left( k+n_{1},l\right) \right) \left(
1-v_{23}\left( l,m\right) \right) & =0  \tag{ii\thinspace a}
\label{Thm5.1(c)(2a)} \\
\left( v_{23}\left( l,m\right) -v_{23}\left( l,m+n_{3}\right) \right) \left(
1-v_{12}\left( k,l\right) \right) & =0,  \tag{ii\thinspace b}
\label{Thm5.1(c)(2b)}
\end{align}
and 
\begin{align}
\left( v_{13}\left( k,m\right) -v_{13}\left( k,m+n_{3}\right) \right) \left(
1-v_{12}\left( k,l\right) \right) & =0  \tag{iii\thinspace a}
\label{Thm5.1(c)(3a)} \\
\left( v_{12}\left( k,l\right) -v_{12}\left( k,l+n_{2}\right) \right) \left(
1-v_{13}\left( k,m\right) \right) & =0. 
\settowidth{\qedskip}{$\displaystyle
  \left( v_{13}\left( k,m\right) -v_{13}\left( k,m+n_{3}\right) \right) \left(
  1-v_{12}\left( k,l\right) \right) =0.
  $}\addtolength{\qedskip}{-\textwidth}\makebox[0pt][l]{\makebox[-0.5\qedskip][r]{\qed}\hss}%
\tag{iii\thinspace b}  \label{Thm5.1(c)(3b)}
\end{align}
\renewcommand{\qed}{}\end{proof}%

\begin{example}
\label{Exa5.2}\emph{Not} all the spectral pairs in three dimensions are
quasi-commutative (although this is true in $d=2$). Take for example the
case (\ref{eq17}) of Section \ref{S2} with 
\begin{equation}
\Lambda =\left\{ 
\begin{pmatrix}
                   k \\
                   \beta \left( k\right) +l \\
                   \gamma \left( k,l\right) +m
                   \end{pmatrix}
                   %
:k,l,m\in \mathbb{Z}\right\}   \label{eq50}
\end{equation}
with $\beta \colon \mathbb{Z}\rightarrow \left[ 0,1\right\rangle $ and $%
\gamma \colon \mathbb{Z}^{2}\rightarrow \left[ 0,1\right\rangle $
arbitrarily given functions. Then the three operators $V_{23}$, $V_{13}$ and 
$V_{12}$ are as follows: 
\begin{align}
&V_{23} =I\text{\quad (the identity operator in the two marked tensor slots),%
}  \tag{i}  \label{Exa5.2(1)} \\
&V_{13}\left( e_{k}^{\left( 1\right) }\otimes e_{\gamma \left( k,l\right)
+m}^{\left( 3\right) }\right)  =e^{i2\pi \beta \left( k\right)
}e_{k}^{\left( 1\right) }\otimes e_{\gamma \left( k,l\right) +m}^{\left(
3\right) }\,,  \tag{ii}  \label{Exa5.2(2)} \\
\intertext{and}%
&V_{12}\left( e_{k}^{\left( 1\right) }\otimes e_{\beta \left( k\right)
+l}^{\left( 2\right) }\right)  =e^{i2\pi \gamma \left( k,l\right)
}e_{k}^{\left( 1\right) }\otimes e_{\beta \left( k\right) +l}^{\left(
2\right) }\,.  \tag{iii}  \label{Exa5.2(3)}
\end{align}
It follows that the three commuting unitary one-parameter groups associated
with $\Lambda $, via Theorem \ref{Thm1.1}, are not quasi-commuting if the
two functions $\beta $ and $\gamma $ in formula (\ref{eq50}) are both
non-constant.
\end{example}

\begin{acknowledgements}
The authors gratefully acknowledge excellent typesetting
and graphics production by Brian Treadway.
\end{acknowledgements}

\bibliographystyle{bftalpha}
\bibliography{jorgen}

\end{document}